\input amstex.tex \documentstyle{amsppt}
 
\baselineskip = 17pt plus 2pt \parindent=20pt \magnification\magstephalf
\rm \hsize15truecm \vsize22truecm 
\voffset=12truemm

\def\b#1{\bold #1}


\define\one{\bold {1}} \define\cons{\operatorname{F }}
\define\algc{\operatorname{A}} 
\define\supp{\operatorname{supp }} \define\con{N^* }
\define\as{\Cal A_{\Cal S}}

\define\hlink{\tilde \Lambda } \define\link{ \Lambda } \define\alink{
\Omega } \define\ahlink{\tilde \Omega }

\define\lk{\operatorname{lk}} \define\linkat#1#2#3{\lk_{#1} (#2; #3)}

\topmatter

\title  Algebraically constructible functions \\ \endtitle

\rightheadtext {Algebraically constructible functions} \author  Clint
McCrory and Adam Parusi\accent19nski \endauthor

\address Department of Mathematics, University of Georgia, Athens, GA
30602, USA \endaddress

\email clint\@math.uga.edu \endemail

\address    D\'epartement de Math\'ematiques, Universit\'e d'Angers, 2 bd.
Lavoisier, 49045 Angers Cedex, France, and School of Mathematics and
Statistics, University of Sydney, Sydney, NSW 2006, Australia \endaddress

\email parus\@tonton.univ-angers.fr, parusinski\_a\@maths.su.oz.au
\endemail

\abstract An algebraic version of Kashiwara and Schapira's calculus of
constructible functions is used to describe local topological properties of
real algebraic sets, including Akbulut and King's numerical conditions for
a stratified set of dimension three to be algebraic. These properties,
which include generalizations of the invariants modulo 4, 8, and 16 of
Coste and Kurdyka, are defined using the link operator on the ring of
constructible functions. \endabstract







\thanks Research partially supported by a University of Sydney Research
Grant. First author also partially supported by NSF grant DMS-9403887.
\endthanks

\endtopmatter

\document

In 1970 Sullivan \cite{Su} proved that if $X$ is a real analytic set and
$x\in X$, then the Euler characteristic of the link of $x$ in $X$ is even.
Ten years later, Benedetti and Ded\`o \cite{BD}, and independently Akbulut
and King \cite{AK1}, proved that Sullivan's condition gives a topological
characterization of real algebraic sets of dimension less than or equal to
two. Using their theory of resolution towers, Akbulut and King introduced a
finite set of local ``characteristic numbers'' of a stratified space $X$ of
dimension three, such that $X$ is homeomorphic to a real algebraic set if
and only if all of these numbers vanish \cite{AK2}.

In 1992 Coste and Kurdyka \cite {CK} proved that if $Y$ is an irreducible
algebraic subset of the algebraic set $X$ and $x\in Y$, then the Euler
characteristic of the link of $Y$ in $X$ at $x$, which is even by
Sullivan's theorem, is generically constant mod 4. They also introduced
invariants mod $2^k$ for chains of $k$ strata, and they showed how to
recover the Akbulut-King numbers from their mod 4 and mod 8 invariants. The
Coste-Kurdyka invariants were generalized and given a simpler description
in \cite{MP} using complexification and monodromy.

We introduce a new approach to the Akbulut-King numbers and their
generalizations which is motivated by the theory of Stiefel-Whitney
homology classes, as was Sullivan's original theorem. We use the ring of
constructible functions on $X$, which has been systematically developed by
Kashiwara and Schapira \cite{KS} \cite{Sch} in the subanalytic setting.
Their calculus of constructible functions includes the fundamental
operations of duality and pushforward, which correspond to standard
operations in sheaf theory.

Our primary object of study is the ring of algebraically constructible
functions on the real algebraic set $X$. We say that the function
$\varphi:X\to\Bbb Z$ and the stratification $\Cal S$ of $X$ are {\it
compatible} if $\varphi$ is constant on each stratum of $\Cal S$. If $X$ is
a complex algebraic set, then $\varphi$ is said to be complex algebraically
constructible if there exists a complex algebraic stratification $\Cal S$
of $X$ which is compatible with $\varphi$. The pushforward of a complex
algebraically constructible function by a complex algebraic map is complex
algebraically constructible.

For real algebraic sets the situation is more complicated. By an algebraic
stratification of the real algebraic set $X$, we mean a stratification
$\Cal S$ of $X$ with strata of the form $Y\setminus Y'$, where $Y$ and $Y'$
are algebraic sets. Thus the strata are not necessarily connected. If $X$
is a real algebraic set and $\varphi:X\to\Bbb Z$, let us say that $\varphi$
is {\it strongly algebraically constructible} if there is an algebraic
stratification $\Cal S$ of $X$ which is compatible with $\varphi$. The
pushforward of a strongly algebraically constructible function by an
algebraic map is not necessarily strongly algebraically constructible.
(Consider for example $\varphi=f_*\one_{\Bbb R}$, where $f:\Bbb R\to \Bbb
R$, $f(x)=x^2$. Then $\varphi(x)=2$ if $x>0$, $\varphi(0)=1$, and
$\varphi(x)=0$ if $x<0$.) On the other hand, we say that $\varphi$ is {\it
semialgebraically constructible} if there is a semialgebraic stratification
of $X$ which is compatible with $\varphi$. The pushforward of a
semialgebraically constructible function by a continuous semialgebraic map
is semialgebraically constructible. But information about the algebraic
structure of $X$ is lost by passing to the ring of semialgebraically
constructible functions.

To solve this dilemma we adopt a definition of algebraic constructibility
which is not solely in terms of compatibility with a stratification. If $X$
is a real algebraic set, we say that $\varphi:X\to\Bbb Z$ is {\it
algebraically constructible} if $\varphi$ is the pushforward of a strongly
algebraically constructible function by an algebraic map. It follows that
the pushforward of
 an algebraically constructible function is algebraically constructible;
however, not every semialgebraically constructible function is
algebraically constructible. (For example let $X=\Bbb R$ and let
$\varphi(x)=1$ if $x\geq 0$, $\varphi(x)=0$ if $x<0$.) We detect the
difference between algebraically constructible functions and
semialgebraically constructible functions by means of the topological {\it
link operator} $\Lambda$ on the ring of semialgebraically constructible
functions. The link operator generalizes the link of a point in a space,
and it is related to the Kashiwara-Schapira duality operator $D$ by
$D\varphi=\varphi-\Lambda\varphi$.

Our main results are the following. Using resolution of singularities, we
prove that if $\varphi$ is an algebraically constructible function then
 $\frac 12 \Lambda\varphi$ is algebraically constructible, and in
particular $\frac 12 \Lambda\varphi$ is integer-valued. We give a new
description of the Akbulut-King numbers in terms of the operator
$\widetilde\Lambda =\frac 12 \Lambda $, and we prove that if $X$ is a
semialgebraic set of dimension less than or equal to three, then $X$ is
homeomorphic to an algebraic set if and only if all of the functions
obtained from $\one_X$ by the arithmetic operations $+,-,*$, together with
the operator $\widetilde\Lambda$, are integer-valued.

We prove the basic properties of (semialgebraically) constructible
functions in section 1. We derive some properties of constructible
functions $\varphi$ which are self-dual ($D\varphi=\varphi$) or
anti-self-dual ($D\varphi=-\varphi$). If $\varphi$ is compatible with a
stratification $\Cal S$ which has only even (resp.~odd) dimensional strata,
then $\varphi$ is self-dual (resp.~anti-self-dual). If $\varphi$ is
self-dual (resp.~anti-self-dual), then the Stiefel-Whitney homology classes
\cite{FM} satisfy $\beta w_i(\varphi)=w_{i-1}(\varphi)$ for $i$ even
(resp.~odd), where $\beta$ is the Bockstein homomorphism.

In section 2 we introduce algebraically constructible functions, and we
give examples of functions which are constructible but not algebraically
constructible, and functions which are algebraically constructible but not
strongly algebraically constructible. We prove that if $\varphi$ is
algebraically constructible then $\widetilde\Lambda(\varphi)$ is
algebraically constructible. Also we show that the specialization of an
algebraically constructible function is algebraically constructible. We
prove that if $\varphi$ is a constructible function on an algebraic set of
dimension $d$, then $2^d\varphi$ is algebraically constructible.

A constructible function $\varphi$ is {\it Euler} if the function
$\widetilde \Lambda(\varphi)$ is integer-valued. By a {\it completely
Euler} function we mean a constructible function $\varphi$ such that all
the functions obtained from $\varphi$ by means of the arithmetic operations
$+,-,*$ and the operator $\widetilde\Lambda$ are integer-valued. In section
3 we analyze such functions in low dimensions. We give computable
conditions to determine whether a constructible function $\varphi$ is
completely Euler, in the case that $\varphi$ has support of dimension less
than or equal to 2, and to determine whether $\one_X$ is completely Euler,
in the case that $X$ has dimension less than or equal to 3.

In section 4 we apply the preceding results to the topology of real
algebraic sets. We give a new proof of our theorem \cite{MP} concerning the
iterates of the relative link operator $\Lambda_YX$ for $Y$ an algebraic
subset of $X$: If $X_1,\dots,X_k$ is an ordered collection of algebraic
subsets of $X$, then $\varphi=\Lambda_{X_1}\cdots\Lambda_{X_k}\one_X$ is
divisible by $2^k$, and if $Y$ is an irreducible algebraic subset of $X$,
then $\varphi$ is generically constant mod $2^{k+1}$ on $Y$. We give a new
description of Akbulut and King's necessary and sufficient conditions for a
compact semialgebraic set $X$ of dimension three to be homeomorphic to an
algebraic set. We prove that $X$ satisfies the Akbulut-King conditions if
and only if $\one_X$ is completely Euler. We give a similar description of
Akbulut and King's conditions for a stratified semialgebraic set to be
homeomorphic to a stratified real algebraic set, by a homeomorphism which
preserves the strata.

In section 5 we introduce {\it Nash constructible} functions, and we show
that a closed semialgebraic set $S$ is symmetric by arcs \cite{Ku} if and
only if $\one_S$ is Nash constructible.

An appendix contains proofs of some elementary foundational results.

For the definitions and properties of real algebraic and semialgebraic sets
and maps, and semialgebraic stratifications, we refer the reader to
\cite{BR}. We will always assume that semialgebraic maps are continuous. By
a real algebraic set we mean the locus of zeros of a finite set of
polynomial functions on $\Bbb R^n$.

\vskip 30pt \head 1.  Constructible functions \endhead 

\medskip Let $X$ be a real algebraic set.
A function $\varphi:X\to \b Z$ is
called ({\it semialgebraically}) {\it constructible}  if it admits a
presentation as a finite sum $$ \varphi = \sum m_i  \one_{X_i}, \tag 1.1 $$
where for each $i$, $X_i$ is a semialgebraic subset of $X$, $\one_{X_i}$ is
the characteristic function of $X_i$, and $m_i$ is an integer.  Denote by
$\cons (X)$ the ring of constructible functions on $X$,  with the usual
operations of addition and multiplication.  The presentation (1.1) is not
unique, but one can always find a presentation with all $X_i$ closed in
$X$.  In what follows, unless otherwise stated, we always assume that the
$X_i$ are closed.  If the support of $\varphi$ is compact, then we may
choose  all $X_i$ compact.  Then  the {\it Euler integral} of $\varphi$ is
defined as $$ \int \varphi = \sum m_i \chi (X_i). $$ By additivity of the
Euler characteristic, the Euler integral does not depend on the
presentation (1.1) of $\varphi$, provided all $X_i$ are compact. Suppose
$Y$ is a semialgebraic subset of $X$ such that the intersection of $Y$ with
the support of $\varphi$ is compact.  Then by $\int_Y \varphi$ we mean the
Euler integral of the restriction of $\varphi$ to $Y$.

Let $f:X\to Y$ be a (continuous) semialgebraic map of real algebraic sets.
If $\psi \in \cons (Y)$, the {\it inverse image}, or {\it pullback}, of
$\psi$ by $f$ is defined by $$ f^*\psi (x) = \psi(f(x)), $$ and $f^*\psi$
is a constructible function on $X$.

Assume that $f:X\to Y$ restricted to the support of $\varphi \in \cons (X)$
is proper.
 Then the {\it direct image}, or {\it pushforward}, $f_*\varphi \in \cons
(Y)$ is given by
the formula $$ f_*\varphi (y) = \int_{f^{-1}(y)} \varphi.$$

Suppose that $X$ is embedded in $\b R^n$.  Then we define the {\it link} of
$\varphi$ as the constructible function on $X$ given by $$ \Lambda \varphi
(x) = \int_{S(x,\varepsilon)} \varphi , $$ where $\varepsilon>0$ is
sufficiently small, and $S(x, \varepsilon)$ denotes the
$\varepsilon$-sphere centered at $x$. The function $\Lambda\varphi$ is
independent of the embedding of $X$ in $\Bbb R^n$. This follows from the
fact that the link of a point in a semialgebraic set is well-defined up to
semialgebraic homeomorphism ({\it cf.}~the Appendix). The duality operator
$D$ on constructible functions, introduced by Kashiwara and Schapira in
\cite {KS},  satisfies $$ D \varphi = \varphi - \Lambda \varphi, $$ which
is equivalent to formula (2.7) of \cite{Sch}.

\medskip \proclaim {1.2. Proposition}

(i) $D(D\varphi) = \varphi$,

(ii) $f_* D = D f_* $,

(iii) $(g\circ f)_* = g_* \circ f_* $. \endproclaim

\medskip \demo {Proof} (i)-(iii)
are proved in \cite {Sch} using the corresponding
operations on constructible sheaves.  For a different proof see the
Appendix below. \qed\enddemo

\proclaim {1.3. Corollary}

(i) $\Lambda \circ \Lambda = 2 \Lambda $,

(ii) $f_* \Lambda = \Lambda f_* $,

(iii) $\int \Lambda \varphi = 0 $. \endproclaim

\demo {Proof} (i)-(ii) are clear.  If the support of $\varphi$ is compact
then (iii) follows from (ii).  Indeed, let $f:X \to P$ be a constant map to
a one point space $P$. Then $$ \int \Lambda \varphi = f_* \Lambda \varphi
(P) = \Lambda f_* \varphi (P) = 0, $$ since the link of any constructible
function on $P$ vanishes. In general we need only the compactness of the
support of $\Lambda \varphi$ for (iii) since this case reduces to the
previous one by (i).   \qed\enddemo

\medskip A constructible function $\varphi$
is called {\it self-dual} if $D \varphi
= \varphi$,  or equivalently $\Lambda \varphi =0$. Similarly, $\varphi$ is
{\it anti-self-dual} if $D \varphi = -\varphi$, or equivalently
$\Omega\varphi=0$, where $\Omega\varphi=\varphi+D\varphi$.

We say that $\varphi \in \cons (X)$ is {\it Euler} if $\Lambda \varphi (x)$
is even for all $x\in X$.  Clearly every self-dual and every anti-self-dual
function is Euler.  On the other hand, every Euler function admits a
canonical decomposition into self-dual and anti-self-dual parts, $$ \varphi
= \hlink \varphi + \ahlink \varphi , \tag 1.4 $$ where $\hlink = \frac 1 2
\Lambda$ and $\ahlink = \frac  1 2 \Omega$.  By (ii) of Proposition 1.2 the
direct image of a function which is Euler (resp.~self-dual, anti-self-dual)
is Euler (resp.~self-dual, anti-self-dual).  Note that whether a
constructible function is Euler it depends only on its reduction modulo 2.
This is no longer true for self-dual or anti-self-dual functions.

We list here some more consequences of Proposition 1.2 which we use in the
sequel: $$ \aligned & D\circ \Lambda = - \Lambda \circ D, \quad D\circ
\Omega = - \Omega \circ D \\ &\hlink \circ \hlink = \hlink, \quad \ahlink
\circ \ahlink = \ahlink ,  \quad \hlink \circ \ahlink = \ahlink \circ
\hlink = 0 . \endaligned \tag 1.5 $$

\medskip Let $\Cal S$ be a
semialgebraic stratification of $X$. We say that $\Cal S$ is {\it locally
trivial} if $X$ as a stratified set can be topologically trivialized
locally along each stratum of $\Cal S$. For instance every Whitney
stratification is locally trivial.  Also a semialgebraic triangulation of
$X$ gives rise to a locally trivial stratification of $X$ by taking open
simplices as strata.  We say that $\varphi\in \cons (X)$ and $\Cal S$ are
{\it compatible} if $\varphi$ is locally constant on strata of $\Cal S$.
For each constructible function $\varphi\in \cons (X)$ there exist a
Whitney stratification of $X$ and a triangulation of $X$ that are
compatible with $\varphi$.

Although Proposition 1.2 is elementary, it carries a nontrivial
information.  For instance, (i) of Proposition 1.2 implies the well-known
fact that the links of points in complex algebraic sets have Euler
characteristic zero.  Actually we can show a more general fact:

\medskip \proclaim {1.6. Proposition}  Let $\varphi$ be a constructible
function on $X$ compatible with a locally trivial stratification $\Cal S$.
Then if all strata of $\Cal S$ are of even (resp.~odd) dimension then
$\varphi$ is self-dual (resp.~anti-self-dual).

In particular if all strata of $\Cal S$ are of odd dimension and the
support of $\varphi$ is compact then $\int \varphi = 0$.  \endproclaim

\medskip \demo{Proof} We show the even-dimensional case.  The proof is by
induction on the dimension of the support $\supp \varphi$.  First note that
the proposition holds generically on $\supp \varphi$. Indeed the geometric
links of a single stratum are odd-dimensional spheres and have zero Euler
characteristic.

Next, by the assumption on local topological triviality,  $\Cal S$ is also
compatible with  $\Lambda \varphi$.  But, by our previous observation, the
dimension of $\supp \Lambda \varphi$ is strictly smaller than the dimension
of $\supp \varphi$.  Hence by inductive hypothesis the statement holds for
$\Lambda \varphi$; that is, $\Lambda \Lambda \varphi = 0$, which by virtue
of (i) of Corollary 1.3 implies $\Lambda \varphi = 0$. This completes the
proof of the even-dimensional case.

The proof in the odd-dimensional case is similar and uses $\Omega$ instead
of $\Lambda$.  The last statement follows from Corollary 1.3 (iii).
\qed\enddemo

In contrast to the direct image, the inverse image does not have good
functorial properties.  In particular, it commutes neither with the duality
operator nor with the link operator.  The following proposition, which we
prove in the Appendix,  shows that the restriction to a generic slice and
the duality operator anticommute.
   
\proclaim {1.7. Proposition}   Let $h:X\to \b R$ be semialgebraic and
let $\varphi\in \cons (X)$.  Let $\varphi_t$ denote the restriction of
$\varphi$ to the fibre $X_t= h^{-1} (t)$.  Then for generic $t\in \b R$ we
have $$ (D \varphi)_t = - D \varphi_t, \quad (\Lambda  \varphi)_t = \Omega
\varphi_t, \quad (\Omega  \varphi)_t = \Lambda \varphi_t. \qed $$
\endproclaim

\medskip Let $f:X\to \b R$ be semialgebraic, and let $x\in X_0 = f^{-1}
(0)$.  Fix a local semialgebraic embedding $(X,x) \subset (\b R^n, 0)$.
Then we define the {\it positive}, resp.~{\it negative}, {\it Milnor fibre}
of $f$ at $x$ by $$ \aligned &   F_f^+(x)  =  B(0, \varepsilon) \cap f^{-1}
(\delta), \\ &   F_f^-(x) =  B(0, \varepsilon )\cap f^{-1} (-\delta) ,
\endaligned $$ where  $B(0, \varepsilon )$ is the ball of radius
$\varepsilon$ centered at $0$ and $0<\delta\ll \varepsilon\ll 1$.

Let $\varphi \in \cons (X)$.  We define the {\it positive} (resp.~{\it
negative}) {\it specialization} of $\varphi$ with respect to $f$ by $$
\aligned &   (\Psi^+_f \varphi) (x) = \int_{F_f^+(x)} \varphi,   \\ &
(\Psi^-_f \varphi) (x) = \int_{F_f^-(x)} \varphi . \endaligned $$ Both
specializations are well-defined, and they are constructible functions
supported in $X_0$.

If $Y$ is a closed semialgebraic subset of $X$, then there exists, at least
locally, a non-negative semialgebraic function $f:X \to \b R$ such that
$Y=f^{-1} (0)$.  For instance, if $X$ is a subset of $\b R^n$ then we may
take $f$ to be the distance to $Y$.  If $x\in Y$, then by the {\it link
along $Y$ at $x$}, denoted $\linkat x Y X$, we mean the positive Milnor
fibre of $f$ at $x$.  If $\varphi$ is a constructible function on $X$, by
the {\it link of $\varphi$
 along $Y$} we mean the positive specialization of $\varphi$ with respect
to $f$, denoted by $\Lambda_Y \varphi $.

The link of $X$ at $x\in X$, which we denote by $\lk (x;X)$, is
well-defined up to semialgebraic homeomorphism, as proven in \cite {CK,
Prop. 1} using \cite {SY}. A similar argument using \cite H shows that the
link of $X$ along $Y$ at $x$ is well-defined up to a semialgebraic
homeomorphism.  A sheaf-theoretic construction of \cite {DS}, see also
\cite {MP, Remark 2},  shows that the   cohomology of $\linkat x Y X$, with
coefficents in any semialgebraically constructible sheaf, is well-defined.
This construction shows that the Euler characteristic of the link is a
topological invariant, which we also prove by elementary means in the
Appendix.

We note also the the Milnor fibres of $f:X\to \b R$ are special cases of
the link construction since $(\Psi^{\pm}_f \varphi) (x) = \linkat x Y
{X_\pm}$, where $X_+ = f^{-1} [0,\infty)$, $X_- = f^{-1} (-\infty,0]$.

\medskip \proclaim {1.8. Proposition}  Let $f:X\to \b R$ be semialgebraic
and continuous. Let $\varphi \in \cons (X)$. Then $\Psi^+_f \varphi +
\Psi^-_f \varphi$ does not depend on $f$ but only on $Y = f^{-1} (0)$ and
equals $$ \Lambda_Y \varphi =  \varphi|_Y - D( (D\varphi)|_Y) = \Lambda (
\varphi|_Y) - \Lambda ( (\Lambda \varphi)|_Y) +  (\Lambda \varphi)|_Y . $$
\endproclaim

\demo {Proof} If one replaces $f$ by $f^2$ one gets $\Psi^+_f \varphi +
\Psi^-_f \varphi = \Lambda_Y \varphi$.  The formula is shown in the
Appendix.  \qed \enddemo

\medskip \proclaim {1.9 Corollary}
Let $f:X\to \b R$ be semialgebraic and
continuous and let $\varphi \in \cons (X)$.  Let $Y = f^{-1} (0)$.  Then $$
\Lambda_Y \circ D = - D \circ \Lambda_Y, \quad \Lambda_Y \circ \Lambda =
\Omega \circ \Lambda_Y, \quad \Omega_Y \circ \Lambda = \Lambda \circ
\Lambda_Y, $$ and similar formulas hold if we replace $\Lambda_Y$ by
 $\Psi^+_f$ or $\Psi^-_f$. \endproclaim

\demo {Proof}  The formulas for $\Lambda_Y$ follow immediately from
Propostion 1.8.  The case of the specializations $\Psi^+_f$ or $\Psi^-_f$
can be reduced to the link $\Lambda_Y$ by replacing $X$ by $X_+ = f^{-1}
[0,\infty)$ or  $X_- = f^{-1} (-\infty,0]$.  \qed \enddemo

Let $X$ be an algebraic subset of the nonsingular real algebraic set $M$.
In \cite {FM} there is defined for each $\varphi \in \cons (X)$ the {\it
conormal cycle} $\con (\varphi)$, which is a Legendrian cycle in $ST^* M$,
the cotangent ray space of $M$. The Euler integral of $\varphi$ can be
computed from the conormal cycle by a generalization of the Gauss-Bonnet
theorem, and the duality operator on constructible functions corresponds to
the action of the antipodal map (multiplication by $-1$ in the fibres) on
the conormal cycle.

Fu and McCrory show that, for each $i=0,1,2,\dots$, there exists a unique
additive natural transformation $w_i$ from $\b Z_2$-valued compactly
supported Euler constructible functions to mod 2 homology,
 such that if $X$ is nonsingular and purely $d$-dimensional, then
$w_i(\one_X)$ is Poincar\'e dual to the classical cohomology
Stiefel-Whitney class $w^{d-i}(X)$.  For $\varphi\in \cons(X)$, the class
$w_i(\varphi)\in H_i(X;\b Z_2)$ is the $i$th {\it Stiefel-Whitney class} of
$\varphi$. Clearly the Stiefel-Whitney classes can be defined for a $\b
Z$-valued constructible function $\varphi$ by first taking the reduction
mod $2$ of $\varphi$.  However by doing that one loses some
information---for instance at the level of $\b Z_2$ coefficients one cannot
distinguish self-dual and anti-self-dual functions.  Instead one may follow
the construction of \cite {FM}, which uses conormal cycles. In particular
one gets the following proposition, which generalizes
 the well-known fact that for  a manifold $M$ of pure dimension $n$,
$w_{n-2k-1}(M)$ can be defined over $\b Z$; this is implied by the fact
that $w_{n-2k-1}(M)$ is the image of $w_{n-2k}(M)$ by the Bockstein
homomorphism $\beta:H_{n-2k}(M;\b Z_2)\to H_{n-2k-1}(M;\b Z_2)$ ({\it cf.}
\cite {HT}).
 
\proclaim {1.10. Proposition}  Let $X$ be a real algebraic set, and let
$\varphi \in \cons (X)$. If $\varphi$ is self-dual, then the odd
dimensional Stiefel-Whitney classes $w_{2k-1}(\varphi)$ are the images of
the even dimensional Stiefel-Whitney classes $w_{2k}(\varphi)$ by the
Bockstein homomorphism.

If $\varphi$ is anti-self-dual then the even dimensional Stiefel-Whitney
classes  are the images of the odd dimensional Stiefel-Whitney classes by
the Bockstein homomorphism. \endproclaim

\demo{Proof} Let $\varphi$ be an Euler constructible function on $X$. Embed
$X$ as an algebraic subset of the smooth $n$-dimensional algebraic set $M$.
Then for $i\geq 0$, the $i$th Stiefel-Whitney class of $\varphi$ is defined
in \cite{FM 4.6} by $$w_i(\varphi)= (\pi_X)_*([\Bbb PN^*
(\varphi)]\frown\gamma^{n-i-1}),$$ where $PN^*(\varphi)$ is the
projectivized conormal cycle in $\Bbb PT^*M$, $\gamma$ is the mod 2 Euler
class of the tautological line bundle on $\Bbb PT^*M$, $\pi: \Bbb PT^*M \to
M$ is the projection, $\pi_X:\pi^{-1}(X)\to X$ is its restriction, and
$[\Bbb PN^*(\varphi)]$ is the mod 2 homology class of $\Bbb PN^*(\varphi)$
in $\pi^{-1}(X)$.

The proof that $\Bbb PN^*(\varphi)$ is a cycle mod 2 \cite{FM, 4.5} hinges
on the fact that $\varphi$ is Euler if and only if $a_*N^*(\varphi)\equiv
N^*(\varphi)\pmod 2$. That proof shows that if
$a_*N^*(\varphi)=N^*(\varphi)$, then $\Bbb PN^*(\varphi)$ lifts to a cycle
with integer coefficients, and hence $\beta[\Bbb PN^*(\varphi)] =0$. Now
$$a_*N^*(\varphi)=(-1)^nN^*(D\varphi),$$ where $a:ST^*M\to ST^*M$ is the
antipodal involution \cite{FM, 3.12}. Suppose that the constructible
function $\varphi$ is self-dual ($D\varphi=\varphi$). If we choose the
embedding $X\subset M$ so that $n=\dim M$ is even, then
$a_*N^*(\varphi)=N^*(\varphi)$, and hence $\beta[\Bbb PN^*(\varphi)]=0$.
Therefore we have $$\align \beta w_i(\varphi) &=\beta(\pi_X)_*([\Bbb
PN^*(\varphi)]\frown\gamma^{n-i-1}) \\ &=(\pi_X)_*\beta([\Bbb
PN^*(\varphi)]\frown\gamma^{n-i-1}) \\ &=(\pi_X)_*([\Bbb
PN^*(\varphi)]\frown\beta\gamma^{n-i-1}) \\ &= \left \{\aligned
(\pi_X)_*([\Bbb PN^*(\varphi)]\frown\gamma^{n-i}) \qquad &n-i-1\ \
\text{odd} \\ (\pi_X)_*([\Bbb PN^*(\varphi)]\frown 0) \ \ \ \ \qquad
&n-i-1\ \ \text{even}\endaligned \right.\\ &=\left\{ \aligned
w_{i-1}(\varphi) \qquad &i\ \ \text{even} \\ 0\ \ \qquad \qquad &i\ \
\text{odd.} \endaligned \right. \endalign$$ Here we use elementary
properties of the Bockstein: $$\align \beta(x\frown u) &=(\beta x\frown
u)+(x\frown\beta u),\\ \beta(u\smile v)&=(\beta u\smile v)+(u\smile \beta
v).\endalign$$ The second equation implies that if $\gamma$ is
1-dimensional, then $\beta(\gamma^k)=\gamma^{k+1}$ if $k$ is odd, and
$\beta(\gamma^k)=0$ if $k$ is even. (Recall that $\beta(\gamma)=\gamma^2$.)

On the other hand, if $\varphi$ is anti-self-dual ($D\varphi =-\varphi$),
we choose $M$ so that $n=\dim M$ is odd. Again
$a_*N^*(\varphi)=N^*(\varphi)$, and the above computation shows that
$$\beta w_i(\varphi)=\left\{\aligned w_{i-1}(\varphi)\qquad &i \ \
\text{odd} \\ 0\ \  \qquad \qquad &i \ \ \text{even,} \endaligned\right.$$
as desired.\qed\enddemo

\vskip 25pt \head 2.  Algebraically constructible functions \endhead

\medskip Let $X$ be a real algebraic set.  In this section we define and
investigate the notion of an algebraically constructible function on $X$.
Of course by a simple analogy to the semialgebraic case one can define an
algebraically constructible function as one admitting a presentation (1.1)
with all $X_i$ algebraic subsets of $X$.  Unfortunately this class of
functions is not preserved by such elementary operations as duality or
direct image by regular mappings.  In this paper we propose to call a
different class of functions algebraically constructible.  Namely, a
function $\varphi: X\to \b Z$ will be called {\it algebraically
constructible} if there exists a finite collection of algebraic sets $Z_i$
and regular proper morphisms ${f_i}:Z_i \to X$ such that $\varphi$ admits a
presentation as a finite sum $$ \varphi = \sum m_i {f_i}_* \one_{Z_i} ,
\tag 2.1 $$ where $m_i$ are integers.   We denote by $\algc (X)$ the ring
of algebraically constructible functions on $X$. The functions which admit
a presentation (1.1) with $X_i$ algebraic will be called {\it strongly
algebraically constructible}.  We note that these two sets of functions on
$X$ coincide if we reduce the coefficients $m_i$ modulo 2.  This follows
easily from the following well-known result.

\proclaim {2.2. Lemma} Let $f:Z\to X$ be a regular morphism of real
algebraic sets, and suppose that $X$ is irreducible.  Then there exists a
proper algebraic subset $Y\subset X$  such that $\chi (f^{-1} (x))$ is
constant modulo 2 on $X\setminus Y$.  \endproclaim

\demo {Proof} See, for instance, \cite {AK2, Proposition 2.3.2}. \qed
\enddemo

The rings $\cons (X)$, $\algc (X)$, and of strongly algebraically
constructible functions are all different if $\dim (X) >0$.  Here are some
examples.

\example {2.3. Examples}  \item {(i)} Let $X=\b R$.  The constructible
function $\varphi \in \cons (\b R)$ is strongly algebraically constructible
if and only if $\varphi$ is generically constant.  On the other hand,
$\varphi\in \algc (\b R)$  if and only if $\varphi$ is Euler or,
equivalently in this case, $\varphi$ is generically constant mod 2. \item
{(ii)} Let $\b P^2 = \b P_{\b R}^2 $ be the real projective plane with
homogeneous coordinates $(x:y:z)$.  Let $f:\b P^2 \to \b R^2$ be given by
$f(x:y:z) = ({\frac {x^2} {x^2 + y^2 + z^2}}, {\frac {y^2} {x^2 + y^2 +
z^2}})$. Then the image of $f$ is the triangle $\Delta$ with vertices
$(0,0)$, $(1,0)$, and $(0,1)$. The pushforward $f_* (\one_{\b P^2})$ is an
algebraically constructible function on $\b R^2$ which equals $4$ inside
$\Delta$, $2$ on its sides, $1$ at the vertices, and $0$ in the complement
of $\Delta$. \item {(iii)} Let $\varphi\in \cons (\b R^2)$ equal twice the
characteristic function of the closed first quadrant. Since $\varphi$ is
even, it is Euler. We show in Remark 2.7 below that $\varphi$ is not
algebraically  constructible. \item {(iv)} Let $f$ be a regular function on
$X$.  Then the sign of $f$ is an algebraically constructible function on
$X$.  Indeed,  let $\widetilde X= \{(x, t) \in X\times \b R\ |\ f(x) =
t^2\}$.  Then $\operatorname{sgn} f = \pi_* \one_{\widetilde X} - \one_X$.
Actually, the signs of regular functions generate the ring $\algc (X)$, as
shown in \cite {PS}.  \endexample

It is clear from the definition that the ring $A(X)$ of algebraically
constructible functions is preserved by the direct image by proper regular
maps.  It is also easy to see that $A(X)$ is preserved by the inverse
image.  We shall show below that $A(X)$ is also preserved by the other
standard operations on constructible functions such as
 duality,  link, and specialization.  To show this we need the following
lemma, which is a consequence of resolution of singularities.
 
\proclaim {2.4. Lemma} Let $\varphi \in \algc (X)$.  Then there exists a
presentation (2.1) of $\varphi$ with all $Z_i$ nonsingular and
pure-dimensional.  \endproclaim

\demo {Proof}  It is sufficient to find such a presentation for $\varphi =
\one_Z$ with $Z$ irreducible.  We proceed by induction on $\dim (Z)$.  By
resolution of singularities there exists a proper regular morphism
$\sigma:\widetilde Z\to Z$ with the following properties: $\widetilde Z$ is
irreducible, nonsingular, and of pure dimension, and there is a proper
algebraic subset $\Sigma \in Z$ such that $\sigma$ induces an isomorphism
between $\widetilde Z \setminus \sigma ^{-1} (\Sigma)$ and $Z\setminus
\Sigma$.  Let $\widetilde \Sigma = \sigma ^{-1} (\Sigma)$.  Then $\dim
(\Sigma) < \dim (Z)$ and $\dim (\widetilde \Sigma) < \dim (Z)$.  Finally $$
\one_Z = \sigma _* \one_{\tilde Z} + ( \one _\Sigma - \sigma_*
\one_{\widetilde \Sigma}) , $$ and the second summand admits the required
presentation by the inductive assumption.  This ends the proof.
\qed\enddemo

The next two results are consequences of Lemma 2.4 and Proposition 1.2.
  
\medskip \proclaim {2.5. Theorem} Let $\varphi \in \algc (X)$.  Then
$\varphi$ is Euler and $\tilde \Lambda \varphi \in\algc (X)$.  Hence
$\tilde \Omega \varphi$ and $D\varphi$ are also algebraically
constructible. \endproclaim

\demo {Proof}  Let $\varphi = \sum m_i {f_i}_* \one_{Z_i}$ be a
presentation given by Lemma 2.4.  Then each $Z_i$ is nonsingular and of
pure dimension, and either $\tilde \Lambda \one_{Z_i} = \one_{Z_i}$ if
$\dim Z_i$ is odd or $\tilde \Lambda \one_{Z_i} = 0$ if $\dim Z_i$ is even.
Hence by Corollary 1.3 (ii), $$ \tilde \Lambda \varphi = \tilde \Lambda
\sum m_i {f_i}_* \one_{Z_i} = \sum m_i {f_i}_* \tilde \Lambda \one_{Z_i} =
{\sum}' m_i {f_i}_* \one_{Z_i}, $$ where the latter sum is only over
odd-dimensional $Z_i$.  \qed  \enddemo

\medskip \proclaim {2.6. Theorem} Let $\varphi$ be an algebraically
constructible function on $X$ and let $f:X\to \b R$ be a regular morphism.
Then  $\Psi^+_f \varphi$ and $\Psi^-_f \varphi$ are algebraically
constructible functions on $X$.  \endproclaim

\demo {Proof} Let $\widetilde X\subset X\times \b R$ be the algebraic set
defined by $\widetilde X = \{ (x,t)\ |\ f(x) = t^2\}$.  Let $\pi:
\widetilde X \to X$ denote the standard projection and let
$\widetilde\varphi = \varphi \circ \pi$, $\widetilde f = f\circ
\pi:\widetilde X \to \b R$.  We identify $\widetilde f^{-1}(0)$ with
$X_0=f^{-1}(0)$.

Take $x\in X_0$.  Then the positive Milnor fibre $F^+_{\widetilde f} (x)$
is the disjoint union of two copies of $F^+_f (x)$ and the negative Milnor
fibre $F^-_{\widetilde f}$ is empty.  Hence, by Proposition 1.8, $$
\Psi_{f}^+  \varphi = \tsize\frac 1 2 (\Psi_{\widetilde f}^+ \widetilde
\varphi  +  \Psi_{\widetilde f}^- \widetilde \varphi) = \tsize\frac 1 2
\Lambda_{X_0} \widetilde \varphi =  \tilde \Lambda (
\widetilde\varphi|_{X_0}) - \tilde \Lambda ( (\Lambda
\widetilde\varphi)|_{X_0}) +  (\tilde \Lambda \widetilde \varphi)|_{X_0}.
$$ The function given by the latter expression is algebraically
constructible by Theorem 2.5.  \qed \enddemo

Note that the first part of Theorem 2.5, that is $\varphi$ is Euler, is
equivalent to Sullivan's observation that the links of points in real
algebraic sets have even Euler characteristic.  On the other hand, the
assertion that $\tilde \Lambda \varphi \in \algc (X)$ is much stronger. It
gives further restrictions for a function to be algebraically constructible
that are of ``greater depth''.  Let us consider the simplest possible
example.  On $X= \b R$ the algebraically constructible functions are
exactly those constructible functions that are  Euler, see Example 2.3 (i).
This is no longer true on $\b R^2$.  To see this let us justify the claim
of Example 2.3 (iv).

\remark {2.7. Remark}  Let $\varphi = 2\cdot\one_Q$, where $Q\subset \b R$
is the closed first quadrant.  Since $\varphi$ is even, it is also Euler.
Let $Y\subset \b R^2$  denote the x-axis.
 Then $(\hlink \varphi)|_Y (x)$ equals $1$ for $x\ge 0$ and   $0$ for
$x<0$.  Consequently  $(\hlink \varphi)|_Y$ is not algebraically
constructible and hence, by Theorem 2.5, neither is $\varphi$. \endremark

\medskip On the other hand,
all constructible functions on $\b R^2$ which are
divisible by $4$ are algebraically constructible (by the following
Theorem).  This gives a clear limit to the depth of information carried by
algebraically constructible functions.   In general we have the following
nontrivial fact.

\proclaim {2.8. Theorem} Let $X$ be an algebraic set of dimension $d$. Then
$$ 2^d \cons (X) \subseteq \algc (X) $$ \endproclaim

\demo {Proof}  Let $S\subset X$ be a semialgebraic subset of $X$.  We show
by induction on $d = \dim X$ that $2^d  \one_S$ is algebraically
constructible.  We suppose that $X$ is irreducible.

Up to a set of dimension $<d$, the set $S$ is a finite union of basic open
semialgebraic sets; that is, sets of the form $$ \{ x\in X\ |\ g_1(x)> 0,
\ldots , g_k (x)>0\}, $$ where the $g_i$'s are polynomials on $\b R^n$
\cite{BCR, 2.7.1}. Since a finite intersection of basic open sets is still
basic, for the inductive step it suffices to consider $S$ basic and open.
Then, by \cite {BCR, Th\'eor\`eme 7.7.8}, there exist polynomials $f_1,
\ldots, f_d$ such that $$ \Cal U=\Cal U (f_1,\ldots, f_d) = \{ x\in X\ |\
f_1(x)> 0, \ldots , f_d (x)>0\} $$ is contained in $S$ and $\dim (S
\setminus \Cal U) <d$.  Hence, by the inductive assumption, it suffices to
show $2^d \one_{\Cal U}\in \algc (X)$.
 Let $\widetilde X\subset X\times \b R^d$ be given by $$\widetilde X= \{(x,
t_1,\ldots t_d) \in X\times \b R^d\ |\ f_1(x) = t_1^2, \ldots, f_d(x)=
t_d^2\}.$$   Let $\pi : \widetilde X\to X$ denote the standard projection.
Then $\widetilde Y = \{ f_1(\pi(\tilde x))= \cdots =f_d(\pi(\tilde x))=0\}$
is an algebraic subset of $\widetilde X$.  In particular, $
\one_{\widetilde X \setminus \widetilde Y}$ is algebraically constructible,
and so is $$ \pi_*  \one_{\widetilde X \setminus \widetilde Y} = 2^d
\one_{\Cal U} $$ as required. \qed \enddemo

\medskip \definition {2.9 Definition}  A constructible function $\varphi$
(respectively a set $\Cal F$ of constructible functions) is {\it completely
Euler} if
 all constructible functions obtained from $\varphi$ (resp.~from the
functions in $\Cal F$) by means of the arithmetic operations $+,-,*$, and
the operator $\hlink$, are  integer valued.   \enddefinition

In particular Theorem 2.5 implies that every algebraically constructible
function is completely Euler.  We shall study some consequences of this
fact in section 4.

The following result is an immediate consequence of Theorem 2.8.  In the
next section we give  an alternative purely topological proof of a slightly
more general statement (Proposition 3.1).

\proclaim {2.10. Corollary} Let $\varphi \in \cons (X)$ be divisible by
$2^{\dim X}$.  Then $\varphi$ is completely Euler. \qed \endproclaim

\vskip 30pt \head 3.  Completely Euler functions \endhead

\medskip In this section
we study the completely Euler functions in low dimensions.
In particular we show how to decide in a systematic way whether a
constructible function is completely Euler.  In the next section we show
that the characterization of completely Euler functions obtained in this
way is equivalent to some combinatorial conditions discovered by Akbulut
and King \cite {AK2}.

If $X$ is a semialgebraic set, we denote by $\Cal A_X$ the ideal of $\cons
(X)$ consisting of all $\varphi$ such that for each positive integer k,
$\dim \supp (\varphi \pmod {2^k}) <k$; that is to say, $\varphi$ is
divisible by $2^k$ in the complement of a subset of dimension $<k$.  If $X$
is an algebraic set then by Theorem 2.8 all functions in $\Cal A_X$ are
algebraically constructible.  In particular they are also completely Euler,
which is also a consequence of  the following.

\proclaim {3.1. Proposition} $\Cal A_X$ is preserved by $\hlink$.
\endproclaim

\demo {Proof} We proceed by induction on $d= \dim X$.  Let $\varphi \in
\Cal A_X$, and let $\Cal S$ be a semialgebraic stratification of $X$
compatible with $\varphi$.  Denote by $X^{d-1}$ the $(d-1)$-skeleton of
$\Cal S$, that is the union of strata of dimension $<d$.  Then $\psi =
(\one_X - \one_{X^{d-1}})\varphi$ is divisible by $2^{d}$ and hence both
$\hlink \psi$ and $\ahlink \psi$ are divisible by $2^{d-1}$.  And either
$\hlink \psi$ for $d$ even, or $\ahlink \psi$ for $d$ odd, has support in
$X^{d-1}$.  In both cases $\hlink \psi = \psi - \ahlink \psi \in \Cal A_X$.

On the other hand $\varphi|_{X^{d-1}} =  (\one_{X^{d-1}})\varphi$ is in
$\Cal A_{X^{d-1}}$ and hence satisfies the inductive assumption. Therefore
$\hlink \varphi = \hlink \psi + \hlink (\varphi|_{X^{d-1}}) \in \Cal A_X$,
as required. \qed \enddemo

Fix a constructible function $\varphi$ (or a set of functions $\Cal F$) on
$X$.   We denote by $\hlink\{\varphi\}$ (or $\hlink \Cal F$, respectively)
the set of functions, which in general may not be integer-valued, obtained
{}from $\varphi$ (resp. $\Cal F$) by means of the arithmetic
operations $+,-,*$,
and the operator $\hlink$.  Let $\Cal S$ be a topologically
trivial semialgebraic stratification compatible with $\varphi$. Let
$\cons_{\Cal S}(X)$ be the subring of $\cons(X)$ consisting of functions
compatible with $\Cal S$, and let $\Cal A_{\Cal S}= \Cal A_{\Cal X} \cap
\cons_{\Cal S}(X)$.  By Proposition 3.1, $\as$ is preserved by $\hlink$ and
hence $\as$ is completely Euler. Consequently, in order to determine
whether $\varphi$ is completely Euler
 we may work in $\cons_{\Cal S}(X)$ modulo $\as$;  that is to say, whether
$\varphi$ is completely Euler is determined by its values mod $2^k$ on
strata of dimension k.  In particular, since $\Cal S$ is finite there are
finitely many conditions to check.
  
We begin with some elementary observations.  First note that whether
$\varphi$ is Euler depends only on the reduction of $\varphi$ modulo 2.
Morever, since all positive powers of $\varphi$ are congruent mod 2, $$
\varphi \equiv \varphi^2 \equiv \varphi ^3 \equiv \cdots \pmod 2, \tag 3.2
$$ and if one of them is Euler so are all the others.
 
Let $\dim \supp \varphi \le 1$.  Then whether $\varphi$ is  completely
Euler is determined by its values modulo 2.  Assume that $\varphi$ is
Euler. Then by (3.2) all the powers of $\varphi$ are also Euler. Moreover,
by dimension assumption, $\ahlink \varphi$ has finite support and so
belongs to $\as$.  Hence $\hlink \{\varphi\}$ modulo $\as$ contains at most
one element, namely the class of $\varphi$.  This shows the following
result.

\medskip \proclaim {3.3 Lemma} If  $\dim \supp \varphi \le 1$ then
$\varphi$ is completely Euler if and only if $\varphi$ is Euler.  \qed
\endproclaim

\medskip Note also
that $\dim \supp \hlink \one_X < \dim X$ (for $\dim X$ even) or
$\dim \supp \ahlink \one_X < \dim X$ (for $\dim X$ odd).    Hence the
following observation will allow  a reduction of dimension.

\medskip \proclaim {3.4 Lemma} $\one_X$ is completely Euler if and only if
$\one_X$ is Euler and $\hlink \one_X$ (or equivalently $\ahlink \one_X$) is
completely Euler.  In particular, if $\dim X \le 2$ then $\one_X$ is
completely Euler if and only if $\one_X$ is Euler. \endproclaim

\medskip \demo {Proof} The first
statement is obvious since multiplication by
$\one_X$ acts trivially on $\hlink \{\one_X\}$.  Suppose that $\dim X \le
2$.  Then $\dim \supp \hlink \one_X \le 1$.  If $\one_X$ is Euler, so is
$\hlink \one_X$, since $\hlink \circ \hlink = \hlink$. So the second
statement follows from Lemma 3.3. \qed \enddemo

On the other hand there exist Euler constructible functions $\varphi$ with
$\dim \supp \varphi = 2$ which are not completely Euler (see Example 3.13
below).  Let us consider this case in detail.  We assume $\varphi$ is Euler
and, as before,
 determine $\hlink \{\varphi\}$ modulo $\as$, in particular modulo 4. The
algebra of powers $\varphi, \varphi^2, \ldots$, modulo 4, is generated
(additively) by $\varphi, \varphi^2, \varphi^3$. By (3.2) all these powers
are Euler.   The supports of $\hlink \varphi ^k$, $k=1,2,3$, are contained
in $X^1$, that is the union of strata of dimension $\le$ 1.  Hence, again
modulo $\as$, $\hlink \{\varphi\}$ is generated additively by the products
of the following functions: $$ \varphi,\ \hlink \varphi,\ \hlink
\varphi^2,\ \hlink \varphi^3 . $$ Moreover, all such products except the
powers of $\varphi$ are supported in $X^1$ and hence it suffices to
consider their values mod 2.  Consequently, by (3.2) only the following
products matter: $$ \varphi^a (\hlink \varphi)^b (\hlink \varphi^2)^c
(\hlink \varphi^3)^d, \tag 3.5 $$ where $a,b,c,d = 0 \text { or } 1$ , and
$b +c +d >0$.  Moreover, $\hlink \varphi$, $\hlink \varphi^2$, $\hlink
\varphi^3$ are automatically Euler.  Thus we have the following result.

\proclaim {3.6 Proposition}  If $\dim \supp \varphi \le 2$,  then $\varphi$
is completely Euler if and only if $\varphi$ is Euler and all 11 functions
supported in the one dimensional set $X^1$ and given by (3.5) with $a,b,c,d
= 0 \text { or }1$, $a+b+c+d \ge 2$, are Euler.  \qed \endproclaim

\medskip Suppose $\dim X \le 3$.
By Lemma 3.4,  Proposition 3.6 applied to $\varphi
= \ahlink \one_X$ gives a criterion for $\one_X$ to be completely Euler.
Since in this case $\hlink \varphi= \hlink \ahlink \one_X = 0$ , whether
$\varphi$ is completely Euler is determined by the products $$ \varphi^a
(\hlink \varphi^2)^b (\hlink \varphi^3)^c \tag 3.7 $$ with $a,b,c = 0 \text
{ or } 1$.  Six of these products,  for $b+c>0$, have support in $X^1$. The
functions $\hlink \varphi^2$ and $\hlink \varphi^3$ are automatically
Euler.  Consequently we have the following.

\proclaim {3.8 Proposition} If $\dim X\le 3$, then $\one_X$  is completely
Euler if and only if $\one_X$ is Euler and  the following functions
supported in $X^1$  are Euler: $$ \varphi (\hlink \varphi^2),\ \varphi
(\hlink \varphi^3),\ (\hlink \varphi^2)(\hlink \varphi^3),\ \varphi (\hlink
\varphi^2) (\hlink \varphi^3) , \tag 3.9 $$ where $\varphi =\ahlink
\one_X$.  \qed\endproclaim
 \bigskip
 
The conditions given by Proposition 3.8 can be expressed in an equivalent
way in terms of characteristic sets.   For every $\delta = (\delta_0,
\delta_1, \delta_2) \in (\b Z_2)^3$ define $$ X_{ \delta} =
  \{ x\in X\ |\  \varphi (x) \equiv \delta_0,\
  \hlink \varphi^2 (x) \equiv \delta_1,\
  \hlink \varphi^3 (x) \equiv \delta_2  \pmod 2 \} $$ Note that the
$X_\delta$ are disjoint, not necessarily closed, and of dimension $\le 1$
if $\delta_1 \ne 0$ or $\delta_2 \ne 0$.  The supports of the functions of
(3.7), considered modulo 2, are unions of the sets $X_\delta$.  In
particular the six functions of (3.7) with $b+c >0$ correspond to the six
set $X_\delta$ with  $\delta_1 \ne 0$ or $\delta_2 \ne 0$: $$ \aligned &
\supp_2 \hlink \varphi^2 = X_{111}\cup X_{110}\cup X_{011}\cup X_{010}, \\
& \supp_2 \hlink \varphi^3 = X_{111}\cup X_{101}\cup X_{011}\cup X_{001},
\\ &\supp_2  \varphi (\hlink \varphi^2) = X_{111} \cup X_{110} , \\ &
\supp_2  \varphi (\hlink \varphi^3) = X_{111} \cup X_{101} , \\ & \supp_2
(\hlink \varphi^2) (\hlink \varphi^3) = X_{111} \cup X_{011} , \\ & \supp_2
\varphi (\hlink \varphi^2) (\hlink \varphi^3) = X_{111}, \endaligned \tag
3.10 $$ where by $\supp_2$ we mean the support modulo 2. Thus Proposition
3.8 can be reformulated as follows.

\medskip \proclaim {3.8$'$ Proposition} If $\dim X\le 3$, then $\one_X$ is
completely Euler if and only if it is Euler and the subsets $X_{111},
X_{101}, X_{011}, X_{110}$ of  $X^1$  are Euler.  \qed\endproclaim \bigskip

\remark {3.11 Remark} If $X$ is Euler then  $\supp_2 \hlink \varphi^2$ and
$\supp_2 \hlink \varphi^3$ are Euler. Therefore we may choose in
Proposition $3.8'$ another family of four characteristic sets $X_\delta$,
provided that if these sets are Euler then all the sets $X_\delta$ are
Euler.  For instance, $X_{111},  X_{101}, X_{001}, X_{110}$ is such a
family, which we use in the next section. \endremark \bigskip

Recall that we have fixed a stratification $\Cal S$ of $X$.  Let $X^i$
denote the $i$-skeleton of $\Cal S$ and suppose, in addition,  that all
skeleta of $\Cal S$ are Euler.  We may apply the above method to obtain a
stratified version of Proposition 3.8 that is a characterisation of those
$\Cal S$ such that the family $\{\one_{X^i}\ |\ i=0,1,2,3\}$ is completely
Euler.

\medskip \proclaim {3.12 Proposition}
Let $\Cal S$ be a locally topologically
trivial stratification of a semialgebraic set $X$, $\dim X\le 3$.  Then the
family $\{\one_{X^i}\ |\ i=0,1,2,3\}$ of characteristic functions of the
skeleta of $\Cal S$ is completely Euler if and only if all $\one_{X^i}$ are
Euler and  one of the following equivalent conditions holds: \item {(i)}
The following 12 functions supported in $X^1$ are Euler: $$ \varphi
\one_{X^1}, \ \ \varphi^a (\hlink \varphi^2)^b (\hlink \varphi^3)^c (\hlink
\one_{X^2})^d , $$ where $\varphi =\ahlink \one_X$, $a,b,c,d = 0$ or $1$,
and we consider only $d=0$, $a+b+c \ge 2$, and $d=1$, $a+b+c >0$.  \item
{(ii)}  The following 12 characteristic sets contained in $X^1$ are Euler:
$$\align X_{ \delta,0} &=  \{ x\in X\ |\ \varphi (x) \equiv \delta_0,\
  \hlink \varphi^2 (x) \equiv \delta_1,\
  \hlink \varphi^3 (x) \equiv \delta_2,\
  \hlink \one_{X^2} \equiv 0 \pmod 2 \} ,\\ \text{for}\ \ \delta &=
(1,1,1),\ (1,0,1),\ (0,1,1),\ (1,1,0),\\  X_{ \delta, 1} &=  \{ x\in X\ |\
\varphi (x) \equiv \delta_0,\
  \hlink \varphi^2 (x) \equiv \delta_1,\
  \hlink \varphi^3 (x) \equiv \delta_2,\
  \hlink \one_{X^2} \equiv 1 \pmod 2 \} ,\\ \text{for}\ \ \delta &\ne
(0,0,0),\ \text{and}\\ X' &= \{x\in X_1\ |\
 \varphi (x) \equiv 1,\
  \hlink \varphi^2 (x) \equiv \hlink \varphi^3 (x) \equiv
  \hlink \one_{X^2} \equiv 0 \pmod 2 \} .\endalign  $$ \endproclaim

\demo {Proof} First note that the family $\{\one_{X^i}\ |\ i=0,1,2,3\}$ is
completely Euler if and only if $\{\varphi, \one_{X^2}, \one_{X_1}\}$ is
completely Euler.  The latter family is supported in $X^2$, so we work
modulo 4.  By repeating the arguments of the proofs of Propositions 3.6 and
3.8, we see that $\{\varphi, \one_{X^2}, \one_{X_1}\}$ is completely Euler
if and only if the functions $$ \varphi^a (\hlink \varphi^2)^b (\hlink
\varphi^3)^c (\hlink \one_{X^2})^d  (\one_{X^1})^e , $$ $a,b,c,d,e = 0$ or
$1$, are Euler.  The supports of $\hlink \varphi^2$, $\hlink \varphi^3$,
and $\hlink \one_{X^2}$ are contained  in $X^1$, so if $ b+c+d>0$ we may
forget the last factor.

Since $\hlink \varphi^2$, $\hlink \varphi^3$, and $\hlink \one_{X^2}$ are
automatically Euler, we are left with exactly the 12 functions of condition
(i).

The equivalence of (i) and (ii) can be shown in exactly the same way as
Proposition $3.8'$. \qed\enddemo

\example{3.13 Example} Let $X$ be Akbulut and King's first published
example of an Euler space which is not homeomorphic to a real algebraic set
\cite{Ki, Example, p.~647}. Recall that $X$ is the suspension of the
algebraic set $Y$ shown in Figure 3 ({\it loc.~cit.}, p.~646). Let $A$ be
the suspension of the figure eight, with suspension points $a$, $a'$; let
$B$ be the suspension of three points, with suspension points $b$, $b'$;
let $C$ be an arc with endpoints $c$, $c'$. The space $Y$ is obtained from
the disjoint union of $A$, $B$, $C$, by identifying $a'$ with $b$, $b'$
with $c$, and $c'$ with $a$. (Note that there is a mistake in the picture
of $Y$ in \cite{BR, p.~181.}) In fact $Y$ is homeomorphic to an algebraic
set in projective 3-space, the union of the umbrella $wx^2=yz^2$ and the
circle $x=0$, $(y-1)^2+z^2=w^2$. The support of $\varphi=\tilde
\Omega\one_X$ is of dimension 2 and $\varphi$ is Euler, but $\varphi$ is
not completely Euler. In fact $\tilde\Lambda(\varphi^2)$ is not Euler,
which is exactly the reason that $X$ is not homeomorphic to an algebraic
set.\endexample

\vskip30pt \head 4.  Topology of real algebraic sets \endhead 

\medskip Let $X$ be a triangulable
topological space such that the one point
compactification of $X$ is also triangulable.  By a theorem of Sullivan
\cite {Su}, a necessary condition for $X$ to be homeomorphic to a real
algebraic set is that $X$ is mod 2 Euler space; that is, the Euler
characteristic of the link of every point of $X$ is even.  By \cite {AK1},
\cite{BD} this condition is also sufficient if $\dim X\le 2$, but this is
no longer true if $\dim X=3$.  In this case necessary and sufficient
topological conditions were given by Akbulut and King \cite {AK2}, and then
reinterpreted by Coste and Kurdyka \cite {C}, \cite {CK}. More restrictions
on the Euler characteristic of links of real algebraic sets were given in
\cite {C},  \cite {CK}, and \cite {MP}.  We show below that all these
conditions are simple consequences of Theorem 2.5.

It will be convenient for us to proceed using the language of semialgebraic
geometry. Alternatively, one could use Euclidean simplicial complexes or
subanalytic sets.

Let $X$ be an algebraic subset of $\b R^n$.  Let $Y\subset X$ be closed and
semialgebraic, and let $Z\subset X$ be semialgebraic.  Choose a nonnegative
continuous semialgebraic function $f:X\to \b R$ defining $Y$; that is, $Y=
f^{-1}(0)$.  Recall from section 1 that by the link $\linkat p Y Z$ of $Y$
in $Z$ at $p\in Y$ we mean the positive Milnor fibre of $f|_Z$ at $p$. Such
a link can be understood as a generalization of the link considered in
\cite {C}, \cite {CK}, which was only defined at generic points of $Y$.  In
particular the Coste-Kurdyka link has the same homotopy type as ours; see
\cite {MP, \S 2.3}  for details.   In what follows we use only the Euler
characteristic of the link, that is, the operator $\link _Y$ introduced in
section 1.  In \cite {C} Michel Coste made important observations on the
behaviour modulo 4, 8, and 16, of the Euler characteristic of links of real
algebraic subsets. These results are special cases of the following general
statement.

\proclaim {4.1. Theorem \cite {MP, Theorem 2}} Let $X_1, \ldots, X_k$ be
algebraic subsets of $X$.  Then $\varphi = \link_{X_1} \ldots \link_{X_k}
\one_X $ is always divisible by $2^k$.  Moreover, let $Y$ be an irreducible
algebraic subset of $X$.  Then there exists a proper algebraic subset
$Y'\subset Y$ such that for all $x,x' \in Y\setminus Y'$ $$ \varphi (x)
\equiv \varphi (x') \pmod {2^{k+1}}   . $$ \endproclaim

\demo {Proof} $\varphi /2^k$ is algebraically constructible---and, in
particular, integer-valued---by Proposition 1.8 and Theorem 2.5.  The
second part of the statement follows from Lemma 2.2. \qed \enddemo

In \cite {C}, Theorem 4.1 was shown only for $k=1, 2, 3$, and under special
assumptions.  In particular,  it was assumed that $X_1\subset \cdots
\subset X_k$ and $\dim X = \dim X_k +1 = \cdots = \dim X_1 + k$. This
dimensional assumption was first dropped in \cite {CK, Theorem $1'$} for
$k=1$.  The proof of Theorem 4.1 presented here is different from the proof
in \cite {MP}, which was based on the relation between complex monodromy
and complex conjugation.

In \cite {C} and \cite {CK} the authors show how to use Theorem 4.1 to
recover Akbulut and King's combinatorial conditions \cite {AK2, 7.1.1}
characterizing real algebraic sets of dimension $\le 3$.   We show below
that it is even more natural to look at these conditions as consequences of
Theorem 2.5.
  
\proclaim {4.2. Theorem} Let $X$ be a semialgebraic subset of $\b R^n$
with $\dim X\le 3$.  Then $X$ satisfies the Akbulut-King conditions if and
only if $\one_X$ is completely Euler.  \endproclaim

Thus the main result of \cite {AK2} can be rephrased as follows: If $\dim
X\leq 3$, $X$ is homeomorphic to a real algebraic set if and only if
$\one_X$ is completely Euler.  In particular, Theorem 2.5 shows the
necessity of the Akbulut-King conditions ({\it cf.}~Remark A.7 of the
Appendix).

To prove Theorem 4.2, we first recall the Akbulut-King conditions, using
the approach of \cite {C}, \cite {CK}.  Then we apply the results of
section 3.

Let $X$ be a semialgebraic subset of $\b R^n$, $\dim X\le 3$.  We suppose
that $X$ is Euler, and we fix a locally trivial semialgebraic
stratification $\Cal S$ of $X$.   Let $C_0(X)$ be the union of the
1-skeleton $X^1$ and those  strata $T$ of dimension 2 such that for $x\in
T$, $$ \chi (\linkat x  T X) = \link_T \one_X (x) \equiv 0   \pmod 4 . $$
Equivalently we may say that we include in $C_0(X)$ those two-dimensional
strata $T$ such that $\ahlink \one_X \equiv 1 \pmod 2$ on $T$.  Let
$\varphi = \ahlink \one_X$.  Then in the complement of $X^1$, $$ \varphi
\equiv \one_{C_0(X)} \pmod 2 , \qquad \varphi^2 \equiv \one_{C_0(X)} \pmod
4 . \tag 4.3 $$

\proclaim {4.4 Lemma} $C_0(X)$ is Euler in the complement of $X^0$.
Moreover, for each stratum $S$ of dimension 1 and $p\in S$, $$ \chi
(\linkat p S {C_0(X)}) = \link_S \one_{C_0(X)} (p) \equiv \alink \varphi^2
(p) \pmod 4 .  $$ \endproclaim

\demo {Proof} Fix a stratum $S$ of dimension 1 and let $p$ be a generic
point of $S$.  Let $N$ be a transverse slice to $S$ at $p$.  Denote by $X'$
a small neighbourhood of $p$ in  $N \cap X$, and set $C' = X' \cap C_0(X)$.
Then by Proposition 1.7, $$ \varphi|_{X'} = \hlink \one_{X'} . $$ This
shows $\varphi|_{X'}$ is Euler near $p$ and hence, by (4.3), so is $\one
_{C'}$.  Hence, again by Proposition 1.7,
 $C_0(X)$ is Euler near $p$.  If we apply the same arguments to $\varphi^2$
we get the last equality of the statement.   \qed \enddemo

Given $S$ and $p\in S$ as above, following \cite {C} and \cite {CK} we
consider the number $$ \Delta_p (S,C_0(X),X)   = \chi (\linkat p S X
\setminus \linkat p S {C_0(X)}) - \chi (\linkat p S X) + \chi (\linkat p S
{C_0(X)}) . $$ Note that, as follows from Lemma 4.5 below, $\Delta_p
(S,C_0(X),X)$ does not depend on $p$ but only on $S$ (actually in the
notation of \cite {C} it equals $- \Delta (S,C_0(X),X)$).  Here we follow
the notation of \cite {MP}, where it  is shown that the number $\Delta_p
(S,C_0(X),X)$ has the following geometric interpretation.
 Let $S$, resp.~$C_0(X)$, be given in a neighbourhood of $p$ as the zero
set of a continuous nonnegative semialgebraic function $f$, resp.~$g$.
Then, following \cite {MP}, we define the {\it iterated link}  $\linkat p
{S,C_0(X)} X$ as the iterated Milnor fibre $$ \linkat p {S, C_0(X)} X =
B(p,\varepsilon)\cap f^{-1}(\delta_1) \cap g^{-1}(\delta_2), $$ where $0<
\delta_2\ll \delta_1\ll \varepsilon$.  As shown in \cite {MP, \S 3.4},
$\Delta_p (S,C_0(X),X)$ is the Euler characteristic of $\linkat p
{S,C_0(X)} X$.  Hence \cite {MP, \S 3.5} shows that $$ \Delta_p
(S,C_0(X),X) = \chi (\linkat p {S, C_0(X)} X) = \link_S (\link_{C_0(X)}
\one_X) (p). $$
  
\medskip \proclaim {4.5 Lemma} $$ \Delta_p (S,C_0(X),X) = \alink (\one_
{C_0(X)} \link \one_X) (p). $$ In particular, $$ \Delta_p (S,C_0(X),X)
\equiv 4 \hlink (\varphi^2 + \varphi^3) (p)  \pmod 8 $$ \endproclaim
  
\demo {Proof} We use again a transverse slice $N$ to $S$ at $p$ and
Proposition 1.7.   Let $X'= N\cap X$, $C'=N\cap C_0(X)$, as before.  If $p$
is a generic point of $S$, for instance $S$ is a stratum of a Whitney
stratification of $X$ near $p$, then $$ \link_S (\link_{C_0(X)} \one_X) (p)
= \link ( \link_{C'} \one_{X'})(p). $$ By Proposition 1.8 the expression
above can be written in terms of the link operator $\link$ and the
characteristic functions of $\one_{C'}$ and $\one_{X'}$, after
simplification $\link (\link_{C'} \one_{X'}) = \link (\one_{C'} \alink
\one_{X'})$.  Hence the first formula of the lemma follows again from
Proposition 1.7, since we have to exchange $\link$ and $\alink$ when taking
the slice.  To show the second formula we use  (4.3): $$ \alink (\one_
{C_0(X)} \link \one_X) \equiv   4 \ahlink (\varphi^2 - \varphi^3) \equiv 4
\hlink (\varphi^2 + \varphi^3) + 4(\varphi^2 + \varphi^3) \pmod 8, $$ which
gives the formula since $\varphi^2 + \varphi^3$ is even.  \qed \enddemo

\demo{Proof of Theorem 4.2} Given a 1-dimesional stratum $S$ of $X$, the
Akbulut-King invariant $$(\varepsilon_0(S), \varepsilon_1(S),
\varepsilon_2(S)) \in (\b Z_2)^3$$ is defined as follows (see \cite {C},
\cite {CK}): $$ \aligned \varepsilon_0 (S) & = \tsize\frac 14 \chi (\linkat
p {S, C_0(X)} X ) \hfil \pmod 2  \\ \varepsilon_0 (S) + \varepsilon_1 (S) +
\varepsilon_2 (S) & = \tsize\frac 12 \chi (\linkat p S X) \hfil \pmod 2  \\
\varepsilon_2 (S) & = \tsize\frac 12 \chi (\linkat p S {C_0(X)})\hfil \pmod
2, \endaligned $$ where $p$ can be any point of $S$.  Given $(a,b,c) \in
(\b Z_2)^3$,  define the characteristic set $\Cal E_{abc}(X) $ as the union
of the 0-skeleton $X^0$ and those one-dimensional strata $S$ such that
$(\varepsilon_0(S), \varepsilon_1(S), \varepsilon_2(S))= (a,b,c)$.

Now for every $x\in X$ we define $$ \aligned \varepsilon_0(x) & = \hlink
(\varphi^2 + \varphi^3)(x) \pmod 2 \\ \varepsilon_1(x) & = \hlink
\varphi^3(x) \pmod 2 \\ \varepsilon_2(x) & = \varphi (x) + \hlink \varphi^2
(x) \pmod 2, \endaligned \tag 4.6 $$ where $\varphi=\ahlink\one_X$. If
$(a,b,c)\ne (0,0,0)$, $(0,0,1)$ then $\{x\in X\ |\
(\varepsilon_0(x),\varepsilon_1(x),\varepsilon_2(x))=(a,b,c)\}$ is of
dimension $\le 1$ and hence by Lemmas 4.4 and 4.5, $$ \Cal E_{abc}(X) =
{X^0} \cup  \{x\in X\ |\  (\varepsilon_0(x), \varepsilon_1(x),
\varepsilon_2(x))= (a,b,c)\} .  $$ In particular, for these $(a,b,c)$ the
set $\Cal E_{abc}(X)$ is independent of the choice of stratification (up to
a finite set, since we may always add some point strata).

Note that $(\varepsilon_0(x), \varepsilon_1(x), \varepsilon_2(x))$ equals
$(0,0,0)$,  resp.~$(0,0,1)$, for  nonsingular points of $X$,
resp.~nonsingular points of $C_0(X)$.  In \cite {AK2}, the characteristic
sets corresponding to $(1,1,1)$, $(0,1,0)$, $(1,0,0)$, $(1,1,0)$ are
denoted by $Z_0(X)$, $Z_1(X)$, $Z_2(X)$, $Z_3(X)$ respectively.  It is
shown in \cite {AK2, 7.1.1} that $X$ is homeomorphic to an algebraic set if
and only if $X$ is Euler and $Z_0(X)$, $Z_1(X)$, $Z_2(X)$, $Z_3(X)$ are
Euler. Now the theorem follows from Proposition $3.8'$ and Remark 3.11
since, in the notation of section 3, $\delta_0(x) = \varepsilon_0 (x) +
\varepsilon_1(x) + \varepsilon_2(x)$,  $\delta_1(x) = \varepsilon_0 (x) +
\varepsilon_1(x)$, and $\delta_2(x) = \varepsilon_1(x)$.  Hence $Z_0(X) =
X_{101}$, $Z_1(X) = X_{111}$, $Z_2(X) = X_{110}$, and $Z_3(X) = X_{001}$.
\qed\enddemo \bigskip
 
In \cite {AK2, Theorem 7.1.2} the authors also give a stratified version
of their characterization of real algebraic sets of dimension $\le 3$ which
involves 12 characteristic sets $Z_i(X)$, $i= 0,\ldots,11$.  Again we show
that these combinatorial conditions  follow from Theorem 2.5 and section 3.

\proclaim {4.7 Theorem} Let $\Cal S$ be a locally topologically trivial
semmialgebraic stratification of the semialgebraic set $X$, with  $\dim
X\le 3$, such that all the skeletons $X^i$ of $\Cal S$ are Euler.  Then the
characteristic sets $Z_i(X)$, $i= 0,\ldots,11$, are Euler if and only if
the family $\{\one_{X^i}\ |\ i=0,1,2,3\}$ is completely Euler. \endproclaim

The proof is similar to that of Theorem 4.2.  We just sketch the main
points.

First we recall briefly the construction of $Z_i(X)$, $i= 0,\ldots,11$,
again following ideas of \cite C and \cite {CK}.  It is important to note
that this time the characteristic sets  will depend
 on the stratification $\Cal S$ of $X$. Let $C_0(X)$ be defined as above
and let $C_1(X)$ be  the union of $X^1$ and the remaining 2-dimensional
strata $T$; that is, those strata $T$ such that $\ahlink \one_X \equiv 0
\pmod 2$ on $T$.  Given a 1-dimensional stratum $S$ we define $$
\varepsilon_3 (S) = \tsize\frac 12 \chi (\linkat p S {C_1(X)})\hfil \pmod
2, $$ where $p$ can be any point of $S$.   The following lemma shows that
$\varepsilon_3(S)$ is well-defined.

\proclaim {4.8 Lemma} The set $C_1(X)$ is Euler in the complement of $X^0$.
Moreover, for each stratum $S$ of dimension 1 and $p\in S$, $$
\varepsilon_2 (S) + \varepsilon_3 (S) = \ahlink (\one_{X^2}) (p) \pmod 2.
$$ \endproclaim

\demo {Proof} The set $C_1(X)$ is Euler because so are $C_0(X)$, $X^2=
C_0(X) \cup C_1(X)$, and $X^1= C_0(X) \cap C_1(X)$.  The last statement of
the lemma follows from $$ \varepsilon_2 (S) + \varepsilon_3 (S) =
\tsize\frac 12 \chi (\linkat p S {X^2}) = \ahlink (\one_{X^2}) (p) \pmod 2
. \qed $$ \enddemo

\demo {Proof of Theorem 4.7} Given $(a,b,c,d)\in (\b Z_2)^4$,  define the
characteristic set $\Cal E_{abcd}(X)$ as the union of $X^0$ and those
1-dimensional strata $S$ such that $(\varepsilon_0(S), \varepsilon_1(S),
\varepsilon_2(S), \varepsilon_3(S)) = (a,b,c,d)$. (We follow here the
notation of \cite {AK2, \S 7.1}.)   It is easy to check that  $$ \Cal
E_{abcd}(X) =  X^0 \cup \{x\in X^1\ |\ \ (\varepsilon_0 (x), \varepsilon_1
(x), \varepsilon_2 (x),\varepsilon_3 (x) = (a,b,c,d)\}  , $$ where
$\varepsilon_0 (x)$, $\varepsilon_1 (x)$, $\varepsilon_2 (x)$ are given by
(4.6), and $\varepsilon_3 (x) = \varepsilon_2(x) + \hlink (\one_{X^2}) (x)
\pmod 2 $.  The characteristic sets $Z_i(X)$, $i= 0,\ldots,11$, are unions
of some of the sets $\Cal E_{abcd}(X)$.  The interested reader may consult
\cite {AK2, \S 7.1} for their definitions.  The important property of the
$Z_i$'s  is that they are Euler if and only if all the sets $\Cal
E_{abcd}(X)$ are Euler, as follows from Lemma 7.1.6 {\it loc.~cit}. On the
other hand, it is easy to see by Proposition 3.12 that all the sets $\Cal
E_{abcd}(X)$ are Euler if and only if the family $\{\one_{X^i}\ |\
i=0,1,2,3\}$ is completely Euler.  This completes the proof.

Note also that Proposition 3.12 explains why we need only 12 conditions out
of 16.  \qed\enddemo

\vskip 25pt \head 5.  Nash constructible functions and Arc-symmetric sets.
\endhead

\medskip We present a variation of
our notion of algebraically constructible
functions.

Let  $X$ be a real algebraic set.  A constructible function $\varphi\in
\cons (X)$ is called {\it Nash constructible} if it admits a presentation
as a finite sum $$ \varphi =\sum m_i {f_i}_* \one_{T_i} , $$ where for each
$i$, $m_i$ is an integer, $T_i$ is a connected component of an algebraic
set $Z_i$, and $f_i : Z_i \to X$ is proper and regular.  By the same
arguments as in Section 2, one shows that the family of Nash constructible
functions is preserved by the inverse image by a regular map, the direct
image by a proper regular map, duality, and $\tilde \Lambda$.  Hence not
all constructible functions are Nash constructible. On the other hand,
there are Nash constructible functions which are not algebraically
constructible. Consider for instance the following classical example.  Let
$X\subset \b R^2$ be the curve defined by  $y^2 = (x-1)x(x+1)$.  Then $X$
is irreducible and nonsingular and consists of two connected components
$X_i$, $i=1,2$.  Moreover, the Zariski closure of either of these
components is $X$ itself.  Hence by Lemma 2.2 the characteristic functions
$\one_{X_i}$ are
 not algebraically constructible, though they are clearly Nash
constructible.

 Note that Lemma 2.2 does not hold any longer if we merely assume that $Z$
is a component of a real algebraic set, so this lemma cannot be applied to
study Nash constructible functions.  Instead one can use the following
general statement.

\proclaim {5.1 Proposition} Let $f:Z\to X$ be a proper analytic mapping of
real analytic spaces, and suppose that $X$ is connected and nonsingular.
Then $\chi (f^{-1}(x))$ is generically constant mod 2; i.e., there exists a
subanalytic subset $Y\subset X$ such that $\dim Y < \dim X$ and, for all
$x,x' \in X\setminus Y$, $$ \chi (f^{-1}(x))\equiv \chi (f^{-1}(x')) \pmod
2. $$ Moreover, if $Z$, $X$ and $f$ are semialgebraic, then $Y$ can be
chosen to be semialgebraic. \endproclaim

\demo {Proof}  By \cite {Su} $Z$ is an Euler space.  Let $\varphi = f_*
\one_Z$.  Then $\varphi$ is a subanalytically constructible
 function in the sense of \cite {KS} and \cite {Sch}.  Since, by {\it
loc.~cit.},
 $f_* D = D f_*$, it follows that $\varphi$ is an Euler function; that is,
 $\link \varphi = f_* \link \one_Z$ attains only even values.  Now the
 proposition follows from the following lemma.
 \enddemo
 
\proclaim {5.2 Lemma} Let $X$ be a connected real analytic manifold and
let $\varphi$ be a subanalytically constructible Euler function on $X$.
Then $\varphi$ is generically constant mod 2. \endproclaim

\demo {Proof} $X$ admits a subanalytic triangulation such that $\varphi$ is
constant on open simplices.  Let $\Delta_1$, $\Delta_2$ be two simplices of
dimension $n=\dim X$ such that $\Delta_{12}=\Delta_1\cap \Delta_2$ is a
simplex of dimension $n-1$.  Let $p$ be a point in the interior of
$\Delta_{12}$ and denote by $a_{1}, a_2, a_{12}$ the values of $\varphi$ on
the interiors of $\Delta_1$, $\Delta_2$, and  $\Delta_{12}$ respectively.
Then by definition of the link operator, $\link \varphi (p) = a_1 + a_2 +
(1 +(-1)^{n}) (a_{12} -a_1-a_2)$.  Thus if $\link \varphi (p)$ is even then
$a_1 \equiv a_2 \pmod 2$. \qed \enddemo

\medskip \definition {5.3 Definition} Let $X$ be a real algebraic set.  A
semialgebraic subset $S$ of $X$ is called {\it arc-symmetric} (or {\it
symmetric by arcs}) if,  for every analytic arc $\gamma : (-1,1) \to X$, if
$\gamma ((-1,0))\subset S$ then $\gamma ((-1,1))\subset S$. \enddefinition

Every arc-symmetric semialgebraic set is closed in $X$.  The notion of
arc-symmetric sets was introduced by Kurdyka;   in many ways these sets
resemble algebraic subsets, but they form a much wider class ({\it
cf.}~\cite{Ku}). It is interesting to note that they can be studied using
the techniques introduced in this paper.

\proclaim {5.4. Proposition} Let $S$ be a closed semialgebraic  subset of
an algebraic set $X$.  Then $\one_S$ is Nash constructible if and only if
$S$ is arc-symmetric.  \endproclaim

\proclaim {5.5. Corollary} Every  arc-symmetric semialgebraic set $S$ is
Euler and $\one_S$ is completely  Euler. \qed \endproclaim

\proclaim {5.6. Corollary} Every arc-symmetric semialgebraic set of
dimension $\le 3$ is homeomorphic to an algebraic set.\qed \endproclaim

\demo {Proof of 5.4}  Let $\one_S=\sum m_i {f_i}_* \one_{T_i}$ be Nash
constructible.  Let $\gamma : (-1,1) \to X$  be an analytic arc in $X$ such
that $\gamma((-1,0))\subset S$. Then by Lemma 5.1, $\chi
(f_i^{-1}(\gamma(t)))$ is generically constant mod 2 on $(-1,1)$.  Hence so
is $\one_S$.  This gives, for $S$ closed, $\gamma ((-1,1))\subset S$, as
required.

Conversely, let $S$ be an arc-symmetric semialgebraic subset of $X$.  We
show by induction on $n = \dim S$ that $\one_S$ is Nash constructible.  We
may assume that $X$ is the smallest algebraic set containing $S$; that is,
$X$ is the Zariski closure of $S$.   Then $\dim S = \dim X$.  Let $\sigma :
\widetilde X \to X$ be a resolution of singularities of $X$.  Then $\sigma$
is an isomorphism over $X\setminus \Sigma$, where $\Sigma $ is an algebraic
subset of $X$, and both $\Sigma$ and $\widetilde \Sigma = \sigma ^{-1}
(\Sigma)$ are of dimension smaller than $n$.  Let $\widetilde X_1, \ldots,
\widetilde X_s$ be the connected components of $\widetilde X$ of dimension
$n$, and let $\widetilde S$ be the union of those $\widetilde X_i$ such
that $$ \sigma (\widetilde X_i) \cap (S\setminus \Sigma) \ne \emptyset . $$
Then, since $S$ is arc-symmetric, by an argument of Kurdyka \cite {Ku,
Th\'eor\`eme 2.6},  $$ \sigma (\widetilde S) \subset S . $$ Hence $$ \one_S
= \one_{S\setminus \Sigma} + \one_{S\cap \Sigma}
 = \sigma_* \one_{\tilde S} - \sigma_* \one_{\tilde S\cap \widetilde
\Sigma}
 + \one_{S\cap \Sigma} . $$
         The first two summands are Nash constructible by definition, and
the latter is Nash constructible by the inductive hypothesis, since clearly
$S\cap \Sigma$ is arc-symmetric. \qed \enddemo

\vskip25pt \head Appendix: Proofs of some properties of constructible
functions \endhead

\medskip In this section we present elementary
proofs of some basic properties of
semialgebraically constructible functions. These proofs use either
stratifications or triangulations of semialgebraic sets (see \cite {\L} for
references).  We believe our arguments are well-known to specialists, and
we do not claim any originality ({\it cf.}~\cite {Sch, Remark 3.5}).

Let $X$ be a closed semialgebraic set and let $x\in X$.  By \cite {CK,
Prop. 1} the link $\lk (x; X)$ is well-defined up to semialgebraic
homeomorphism.  The Euler characteristic of the link is a topological
invariant of the germ $(X,x)$.  Indeed, $X$ is locally  contractible and
hence $$ \chi (\lk (x; X)) = 1 - \chi(X,X \setminus \{x\}). \tag A.1 $$

Similarly let $Y$ be a compact semialgebraic subset of $X$.  Then the
quotient space $X/Y$ has a natural structure as a semialgebraic set and we
may define the {\it link of $Y$ in $X$}, by $\lk (Y;X) = \lk (*; {X/Y})$,
 where $*$ denotes the class of $Y$ in $X/Y$.  If $Y\subset \b R^n$ is
compact and semialgebraic, and not necessarily contained in $X$, then by
$\lk (Y;X)$ we mean $\lk (Y\cap X;X)$.  By the above, $\chi (\lk (Y;X))$ is
also a topological invariant of the pair $(X,Y)$; this also follows from
the following corollary of \cite {MP, Lemma 1}: $$ \chi (\lk (Y;X)) = \chi
(Y\cap X) + \chi (X \setminus Y) - \chi (X) =\chi (Y\cap X) - \chi (X,X
\setminus Y). \tag A.2 $$

Let $f:Z\to X$ be a proper semialgebraic map and let $Y$ be a compact
subset of $X$.  In what follows we often use the following consequence of
the definition of the link: $$ f^{-1} (\lk (Y;X)) = \lk (f^{-1}(Y); Z). $$

\medskip Recall that
the additivity of Euler characteristic, $$ \chi (X\cup X') =
\chi (X) + \chi (X') - \chi (X\cap X'), \tag A.3 $$ allowed us to define in
section 1 the Euler integral of the semialgebraically constructible
function $\varphi \in \cons (X)$, provided $\varphi$ has compact support.
Here are some other elementary consequences of (A.3).

\bigskip \proclaim {A.4 Lemma} Let $X$ and $Y$ be closed semialgebraic
subsets of $\b R^n$ and suppose that $Y$ is compact.  Then $$ \chi (\lk (Y;
X)) = \int_Y \link \one_X  . \tag A.4.1 $$ Let $f:Z\to X$ be a
semialgebraic map,
 and let $\varphi \in \cons (Z)$ have compact support.  Then
 $$
 \int f_* \varphi = \int \varphi.
 \tag A.4.2
 $$ Let $x\in X$ and let $\varphi\in \cons (X)$. Then for $\epsilon>0$
sufficiently small, $$ \int_{B_\varepsilon} \varphi = \varphi (x), \tag
A.4.3 $$ where $B_\varepsilon$ is the closed ball of radius $\epsilon$
centered at $x$. \endproclaim

\medskip \demo {Proof}
 The right hand side of (A.2) is additive with respect to $X$ for $Y$ fixed
 and additive with respect to $Y$ for $X$ fixed.  By (A.2) so is $\chi (\lk
(Y;X))$.
 On the other hand, the right hand side of (A.4.1) is also additive with
respect
 to both $X$ and $Y$.   Hence, by the triangulability of semialgebraic
sets,
 it suffices to verify (A.4.1) for $X$ and $Y$ simplices such that
 $X\cap Y$ is
 their common face.  In this case the verification is straightforward.
 This shows (A.4.1).
 
 To show (A.4.2) we may assume that $Z$ is compact and $\varphi = \one_Z$.
 We may assume also that, up to a semialgebraic
 homeomorphism,  $X$ is equal to a simplex $\Delta$ and
 $f$ is topologically trivial with
 fibre $F$ over the interior of $\Delta$.  Let $Y= \partial \Delta$, and
let
 $Z'= f^{-1} (Y)$.  Since
 $f$ is topologically trivial over $\lk (Y; X)$,
 $$
 \chi (\lk (Z'; Z))  = \chi (F) \chi (\lk ( Y; X)) ,
 $$
 which, by (A.2) and the inductive assumption on $\dim X$, gives
 $$
 \aligned
 \int \varphi &= \chi (Z) \\
 &=  \chi (Z') + \chi (Z\setminus Z') -  \chi (\lk (Z', Z)) \\
 &= \int_Y f_* \varphi + \chi (F) (\chi (X\setminus Y)
    -  \chi (\lk (Y; X))) \\
 &= \int_Y f_* \varphi + \chi (F) (\chi (X) -  \chi (Y)) \\
 &= \int \one_Y f_* \varphi + \int \one_{X\setminus Y} f_* \varphi \\
 &= \int f_* \varphi,
 \endaligned
 $$
 as required.  This shows (A.4.2).
 
It suffices to show (A.4.3) for $\varphi = \one_X$ and in this case
(A.4.3) follows from the local contractibilty of $X$.   \qed \enddemo

\medskip \demo {Proof of Proposition 1.2} To show (i) we note that $D$ is
additive, {\it i.e.} $D(\varphi +\psi) = D\varphi + D\psi$.  Therefore it
suffices to verify (i) for $\varphi = \one_\Delta$, where $\Delta$ is a
simplex.  In this case the verification is straightforward.

Let $f:X\to X'$ be proper semialgebraic, $x\in X'$, and $Y = f^{-1} (x)$.
Then by (A.4.1),  $$ \aligned
 f_* (\link \one_X) (x) &= \int_Y \link \one_X = \chi (\lk (Y; X))
 =\chi ( f^{-1} (\lk (x; X')))\\ &= \int_{\lk (x,X')} f_* \one_X = \link
(f_* \one_X) (x) . \endaligned $$ Hence $f_* \link = \link f_*$, which
implies (ii) of Propostion 1.2.

Statement (iii) follows from (A.4.2).  Indeed, let $Z @>f>> X @ >g>> Y$ and
let $y\in Y$.  Then by (A.4.2), $$ (g\circ f)_* \varphi(y) = \int_{(g\circ
f)^{-1}(y)} \varphi = \int_{g^{-1}(y)} f_* \varphi = g_*(f_* \varphi) (y).
\qed $$ \enddemo

\demo {Proof of Proposition 1.7} Let $h:X\to \b R$ be semialgebraic.  Then
$h$ is locally topologically trivial over the complement of a finite subset
of $\b R$.  By \cite {H} we may assume that this trivialization is
semialgebraic.  In particular for $x\in X$ in a generic fibre $X_t$ of $h$,
the link $\lk (x;X)$ is the suspension of $\lk (x;X_t)$.  This gives $$
\chi (\lk (x;X)) = 2 - \chi (\lk (x;X_t)).  \tag A.5 $$  This, in
particular, shows Propostion 1.7.  Note also that, by virtue of (A.1), we
may use any topological trivialization of $h$ (not necessaily
semialgebraic) to establish (A.5).  \qed \enddemo

Suppose, in addition, that $h:X\to \b R$ is proper and let $c_0<c_1$ be
generic values of $h$.  Let $X_{c_0,c_1} = h^{-1}[c_0,c_1]$.  Then by (A.5)
$$ \link \one_{X_{c_0,c_1}} = (\link \one_X)|_{X_{c_0,c_1}}  +
(\one_{X_{c_0}} - (\link \one_X)|_{X_{c_0}}) +
 (\one_{X_{c_1}} - (\link \one_X)|_{X_{c_1}}) . \tag A.6 $$

\medskip \demo {Proof of Proposition 1.8} Let $Y\subset X\subset \b R^n$
and fix $x\in X$. Let $\varphi \in \cons (X)$.  Let $B_\varepsilon$ denote
a small closed ball centered at $x$ and $S_\varepsilon  = \partial
B_\varepsilon$.    By (A.6)  $$ \Lambda (\varphi|_{B_\varepsilon}) =
(\Lambda \varphi)|_{B_\varepsilon} - (\Lambda \varphi)|_{S_\varepsilon} +
\varphi|_{S_\varepsilon},    $$ and hence by (A.4.1) and (A.4.3) $$
\aligned \Lambda_Y \varphi (x) &= \int_{Y\cap B_\varepsilon} \Lambda
(\varphi|_{B_\varepsilon}) \\ &= \int_{Y\cap B_\varepsilon}
[(\Lambda\varphi)|_{B_\varepsilon} - (\Lambda \varphi)|_{S_\varepsilon} +
\varphi|_{S_\varepsilon} ]\\ &= \Lambda\varphi (x)  - \Lambda((\Lambda
\varphi)|_Y) (x)  + \Lambda(\varphi|_Y)(x) \endaligned $$ if $x\in Y$, and
$\Lambda_Y \varphi (x)=0$ otherwise.  This shows Proposition 1.8.
\qed\enddemo

\bigskip \remark {A.7 Remark (Topological invariance of the link operator
and the Euler integral)} Let $h:X'\to X$ be a homeomorphism (not
necessarily semialgebraic) of semialgebraic sets.  Let $\varphi \in \cons
(X)$ be such that $\varphi'= \varphi \circ h \in \cons (X')$.  Let
$Y\subset X$ be a compact semialgebraic subset such that $Y'=h^{-1}(Y)$ is
also semialgebraic.  Then $$ (\Lambda \varphi)\circ h = \Lambda (\varphi'),
\quad \int_Y \varphi = \int_{Y'} \varphi' . $$  Indeed, it suffices to show
that there exist closed semialgebraic sets $X_i\subset X$ such that
$h^{-1}(X_i)$ are semialgebraic subsets of $X'$ and $$ \varphi = \sum m_i
\one_{X_i} . $$  Here is a canonical construction of such sets $X_i$. First
we note that $\varphi$ is semialgebraically constructible if and only if
all the sets $\varphi^{-1}(m)$, $m\in \b Z$, are semialgebraic and all but
finitely many of them are empty.  Let $\varphi_m =
\varphi|_{\varphi^{-1}(m)}$.  Then clearly $$ \varphi =\sum_m \varphi_m .
\tag A.8 $$  Let $Y={\varphi^{-1}(m)}$.  Then $\one_Y$ can be canonically
decomposed $$ \one_Y = \one_{F_1} -  \one_{F_2} + \one_{F_3}- \cdots \pm
\one_{F_d} , \tag A.9 $$ where $F_1\supset F_2 \supset \cdots \supset F_d$
are closed semialgebraic in $\bar Y$.  The sequence of $F_i$'s is
constructed recursively as follows ({\it cf.} \cite {K, \S 12}):  $Y_0 =
Y$, $F_i = \bar Y_{i-1}$, $Y_i= F_i \setminus Y_{i-1}$, $i=1,2,\ldots$.
Clearly all the $F_i$ are closed and semialgebraic, and the sequence
terminates since $\dim F_i < \dim F_{i-1}$.  Now (A.8), together with (A.9)
applied to each set ${\varphi^{-1}(m)}$,  gives the required canonical
decomposition of $\varphi$. \endremark


\bigskip \Refs \widestnumber\key{AGV}

\ref \key AK1 \by     S. Akbulut, H. King \paper The topology of real
algebraic sets    \jour L'Enseignement Math. \vol 29 \yr 1983 \pages
221--261 \endref

\ref \key AK2 \by     S. Akbulut, H. King \book    Topology of Real
Algebraic Sets    \bookinfo  MSRI Publ., vol. 25 \publ    Springer-Verlag
\publaddr New York \yr 1992 \endref

\ref \key BD \by    R. Benedetti, M. Ded\` o \paper The topology of
two-dimensional real algebraic varieties    \jour Annali Math. Pura Appl.
\vol 127 \yr 1981 \pages 141--171 \endref

\ref \key BR \by        R. Benedetti, J-J. Risler \book     Real Algebraic
and Semi-Algebraic Sets \publ      Hermann  \publaddr Paris \yr 1990
\endref

\ref \key BCR \by     J. Bochnak, M. Coste, M-F. Roy \book G\'eom\'etrie
Alg\'ebrique R\'eelle  \publ    Springer-Verlag  \publaddr Berlin \yr 1987
\endref

\ref \key C \by       M. Coste \paper    Sous-ensembles alg\'ebriques
r\'eels de codimension 2      \inbook   Real Analytic and Algebraic
Geometry \bookinfo Lecture Notes in Math.   \publ    Springer-Verlag \yr
1990   \vol      1420 \pages    111--120 \endref

\ref \key CK \by       M. Coste, K. Kurdyka \paper    On the link of a
stratum in a real algebraic set \jour     Topology \vol      31 \yr 1992
\pages 323--336 \endref

\ref \key DS \by        A. Durfee, M. Saito \paper    Mixed Hodge
structures on the intersection cohomology of links \jour     Compositio
Math. \vol   76 \yr 1990  \pages    49--67  \endref

\ref \key FM \by J. H. G. Fu, C. McCrory \paper Stiefel-Whitney classes and
the conormal cycle of a singular variety \jour Trans. Amer. Math. Soc.
\toappear \endref

\ref \key H \by   R. Hardt \paper Semi-algebraic local triviality in
semi-algebraic mappings \jour Amer. Jour. Math. \vol 102 \yr 1978 \pages
291--302 \endref

\ref \key HT \by S. Halperin, D. Toledo \paper Stiefel-Whitney homology
classes \jour Annals of Math. \vol 96 \yr 1972 \pages 511--525 \endref

\ref \key K \by K. Kuratowski \book Topologie \publ Polskie Towarzystwo
Matematyczne \publaddr Warszawa \yr 1952 \vol I \endref

\ref \key Ki \by H. King \paper The topology of real algebraic sets \inbook
Proc. Symp. Pure Math. \vol 40, Part I  \yr 1983 \pages 641--654 \endref

\ref \key Ku \by K. Kurdyka \paper Ensembles semi-alg\'ebriques
sym\'etriques par arcs \jour Math. Ann. \vol 282 \yr 1988 \pages 445-462
\endref

\ref \key KS \by M. Kashiwara, P. Schapira \book Sheaves on Manifolds \publ
Springer-Verlag \publaddr Berlin \yr 1990 \endref


\ref \key \L \by   S. {\L}ojasiewicz \paper Sur la g\'eom\'etrie semi- et
sous-analytique  \jour Ann. Inst. Fourier, Grenoble \vol 43 \issue 5 \yr
1993 \pages 1575--1595 \endref

\ref \key MP \by C. McCrory, A. Parusi\'nski \paper Complex monodromy and
the topology of real algebraic sets \jour Compositio Math. \toappear
\endref

\ref \key PS \by A. Parusi\'nski, Z. Szafraniec \toappear \endref

\ref \key Sch \by       P. Schapira \paper     Operations on constructible
functions  \jour     J. Pure Appl. Algebra \vol     72 \yr 1991 \pages
83--93 \endref

\ref \key SY \by       M. Shiota, M. Yokoi \paper     Triangulations of
subanalytic sets and locally subanalytic manifolds  \jour     Trans. Amer.
Math. Soc. \vol    286  \yr 1984 \pages 727--750 \endref

\ref \key Su \by        D. Sullivan \paper    Combinatorial invariants of
analytic spaces \inbook  Proc. Liverpool Singularities Symposium I
\bookinfo Lecture Notes in Math. \yr         1971 \vol 192 \publ
Springer-Verlag \pages 165--169 \endref

\endRefs
\enddocument